\documentclass[11pt,a4paper]{amsart}
\usepackage[T1]{fontenc}
\usepackage{lmodern,amsmath,amssymb,amsthm,mathtools,array,microtype}
\usepackage[a4paper,hmargin=29mm,top=25mm,bottom=26mm,headheight=14pt,headsep=13pt,footskip=24pt]{geometry}
\usepackage{amscd,fancyhdr,xcolor,needspace,hyperref}
\definecolor{ink}{HTML}{25465D}
\hypersetup{colorlinks=true,linkcolor=ink,citecolor=ink,urlcolor=ink,
pdftitle={}}
\ifdefined\pdfgentounicode\input{glyphtounicode}\pdfgentounicode=1\fi
\fancypagestyle{firstpage}{\fancyhf{}\fancyfoot[C]{\small\thepage}}
\newtheorem{mainthm}{Theorem}

\numberwithin{equation}{section}
\newtheorem{proposition}{Proposition}[section]
\newtheorem{lemma}[proposition]{Lemma}
\newtheorem{corollary}[proposition]{Corollary}
\theoremstyle{definition}
\newtheorem{definition}[proposition]{Definition}
\newtheorem{example}[proposition]{Example}
\theoremstyle{remark}
\newtheorem{remark}[proposition]{Remark}
\newcommand{\K}{\mathcal K}\newcommand{\Zs}{\mathcal Z}
\newcommand{\Z}{\mathbb Z}\newcommand{\Q}{\mathbb Q}\newcommand{\T}{\mathbb T}
\newcommand{\Dsk}{\overline{\mathbb D}}
\newcommand{\Cu}{\operatorname{Cu}}\newcommand{\RR}{\operatorname{rr}}\newcommand{\sr}{\operatorname{sr}}
\newcommand{\rank}{\operatorname{rank}}\newcommand{\im}{\operatorname{Im}}\newcommand{\dist}{\operatorname{dist}}
\newcommand{\diag}{\operatorname{diag}}\newcommand{\Hom}{\operatorname{Hom}}\newcommand{\Ext}{\operatorname{Ext}}
\newcommand{\doi}[2]{\href{https://doi.org/#1}{#2}}
\newcommand{\Tr}{\operatorname{Tr}}

\newcommand{\directlim}{\varinjlim}

\begin{document}
\title[Approximate subhomogeneity in extensions]{Approximate subhomogeneity in extensions of $C^*$-algebras of real rank zero}
\author{Zhichao Liu}
\address{School of Mathematical Sciences, Dalian University of Technology,
Dalian 116024, China}
\email{zhichaoliu@dlut.edu.cn}
\author{Tao Zhang}
\address{School of Mathematical Sciences, Dalian University of Technology,
Dalian 116024, China}
\email{taozhang@mail.dlut.edu.cn}
\date{}

\begin{abstract}
We construct two unital, $\mathcal Z$-stable $C^*$-algebras
of real rank zero with a common simple AT quotient and stable
simple AF ideals. The extensions have the same six-term exact
sequence of abelian groups, with corresponding unit classes,
but different orders on the ideal $K_0$-groups. The first algebra
is an inductive limit of two-dimensional noncommutative CW
complexes and fails cancellation of projections. The second is strongly quasidiagonal, but is not locally
subhomogeneous and its Cuntz semigroup fails near unperforation. 
We also prove that the Cuntz semigroup of every separable,
locally subhomogeneous, $\mathcal Z$-stable $C^*$-algebra is
nearly unperforated.
\end{abstract}

\maketitle
\thispagestyle{firstpage}

\section{Introduction}

 Extension theory links the internal structure of a \(C^*\)-algebra to the interplay between a closed ideal and its quotient. The theory originated in the landmark work of Brown, Douglas, and Fillmore on extensions by the compact operators \cite{BDF}. For a general extension, the associated six-term exact sequence in \(K\)-theory plays a key role through its boundary maps. A fundamental question is how these boundary maps govern the approximation and comparison properties of the middle algebra.

It is a classical result of Elliott \cite[Corollary~3.3]{Ell76} and Brown \cite{BrownAF} that the class of AF algebras is closed under extensions. Real rank zero provides a broader setting in which projections still play a central role, yet boundary obstructions now arise. In this setting, Lin and R{\o}rdam proved that an extension of \(C^*\)-algebras of real rank zero and stable rank one again has stable rank one if and only if the index map vanishes \cite{LR}. Moreover, Brown and Dadarlat established that for any separable AH algebra of real rank zero, every extension arising from a closed ideal has vanishing boundary maps \cite[Corollary~12]{LBD}. Consequently, a nonzero index map precludes an AH presentation, while leaving open the possibility of approximation by broader classes of subhomogeneous algebras.

In the simple setting, Phillips proved that unital inductive limits of recursive subhomogeneous algebras with bounded topological dimension have stable rank one, and hence cancellation \cite{PhillipsLim}. Here we investigate how this conclusion changes when simplicity is relaxed to having a unique nonzero proper ideal. Building on the diagonal constructions for \(1\)-NCCW complexes developed in \cite{AEGL}, we employ a two-dimensional attachment to preserve a nonzero index class in the inductive limit. Although the resulting extension has a simple AF ideal and a simple AT quotient, this index class obstructs cancellation in the middle algebra.

\begin{mainthm}\label{thm:A}
There exists a unital \(\Zs\)-stable ASH algebra of real rank zero which is an inductive limit of \(2\)-NCCW complexes and fails cancellation of projections.
\end{mainthm}
The construction also decouples the algebraic six-term exact sequence from the order structure governing projection comparison. Fixing the quotient, the \(K_0\)-group \(\Q^2\) of the ideal, and the index map, we vary only the positive cone on \(\Q^2\). Consequently, the two resulting extensions share the same algebraic six-term exact sequence and unit class, yet exhibit different comparison behavior.

This phenomenon is directly reflected in the Cuntz semigroup, which extends comparison from projections to positive elements \cite{Cun78,CEI}. While absorption of the Jiang--Su algebra \(\Zs\) implies almost unperforation \cite[Theorem~4.5]{Ror04}, the stronger property of near unperforation means that comparison at all sufficiently large multiplicities already implies comparison at multiplicity one. Antoine, Perera, and Thiel established near unperforation for \(\Zs\)-stable \(C^*\)-algebras that are either simple or of real rank zero and stable rank one, and conjectured that it holds for all \(\Zs\)-stable \(C^*\)-algebras \cite[Corollary~5.6.15 and Conjecture~5.6.18]{APT}. By modifying the positive cone as above, we settle this conjecture in the negative, obtaining a counterexample in which near unperforation already fails for projections and positive elements.
\begin{mainthm}\label{thm:B}
There exists a unital, separable, nuclear, stably finite, \(\Zs\)-stable, real rank zero \(C^*\)-algebra \(\widetilde E\) such that \(\Cu(\widetilde E)\) fails near unperforation. Moreover, there exist two projections \(p,q\in\widetilde E\) with
\[
             [p]_{\mathrm{MvN}}\nleq[q]_{\mathrm{MvN}},\qquad
             n[p]_{\mathrm{MvN}}\le n[q]_{\mathrm{MvN}}\quad(n\ge2).
\]
\end{mainthm}

This example raises a further approximation question. By Hadwin's results \cite{Had}, every separable locally subhomogeneous \(C^*\)-algebra is strongly quasidiagonal, and Brown and Dadarlat used the failure of this property to preclude ASH presentations for certain extensions of real rank zero \cite{LBD}. For the extensions considered here, however, stable finiteness already implies strong quasidiagonality. To obtain a new obstruction, we show that if the middle algebra is locally subhomogeneous, then every trace on the AF ideal vanishes on the image of the index map. For the modified positive cone, the corresponding ideal admits a trace that does not vanish on this image, yielding the following result.

\begin{mainthm}\label{thm:C}
There exists a separable, unital, nuclear, stably finite,
\(\Zs\)-stable, real rank zero \(C^*\)-algebra 
which is strongly quasidiagonal but is not locally subhomogeneous.
\end{mainthm}

The same \(C^*\)-algebra yields both Theorems~\ref{thm:B} and~\ref{thm:C}. Our final result explains the connection between these two phenomena: near unperforation holds for every \(\Zs\)-stable locally subhomogeneous \(C^*\)-algebra.

\begin{mainthm}\label{thm:D}
Let \(A\) be a separable, locally subhomogeneous,
\(\Zs\)-stable \(C^*\)-algebra. Then \(\Cu(A)\) is nearly
unperforated.
\end{mainthm}

The paper is organized as follows. In Section~2, we collect the necessary preliminaries. In Section~\ref{sec:disk}, we construct the extension in Theorem~\ref{thm:A} and compute its ordered \(K\)-theory. In Section~\ref{sec:strict}, we establish a criterion for projection comparison and prove Theorem~\ref{thm:B} by modifying the positive cone on the \(K_0\)-group of the ideal. Finally, in Section~\ref{sec:local}, we prove Theorem~\ref{thm:C} via a trace obstruction to local subhomogeneity, and we prove Theorem~\ref{thm:D} by comparing positive elements after tensoring with \(\Zs\).

\section{Preliminaries}

Denote by \(\K\) the compact operators on a separable infinite-dimensional Hilbert space, and denote by  \(a^{\oplus n}\) the direct sum of \(n\) copies of \(a\). If $A$ is nonunital, we use \(A^\sim\) for its minimal unitization and \(\epsilon\colon A^\sim\to\mathbb C\) for the scalar quotient. For unital \(A\), set \(A^\sim=A\).

\Needspace{9\baselineskip}
\begin{definition}[\cite{Bla}]
For projections \(p,q\in A\otimes\K\), write \(p\sim_{\mathrm{MvN}} q\) if \(v^*v=p\) and \(vv^*=q\) for a partial isometry \(v\), and \(p\precsim_{\mathrm{MvN}} q\) if \(v^*v=p\) and \(vv^*\le q\). Write \([p]_{\mathrm{MvN}}\) for the equivalence class of \(p\). The monoid \(V(A)\) of these classes, with addition induced by direct sum, has algebraic preorder
\[
 [p]_{\mathrm{MvN}}\le[q]_{\mathrm{MvN}}\quad\Longleftrightarrow\quad p\precsim_{\mathrm{MvN}} q.
\]
The algebra \(A\) has \emph{cancellation of projections} if \(p\oplus r\sim_{\mathrm{MvN}} q\oplus r\) implies \(p\sim_{\mathrm{MvN}} q\). It is \emph{stably finite} if no unit of \(M_n(A^\sim)\) is equivalent to a proper subprojection of itself.
\end{definition}

\begin{definition}[\cite{Bla}]
For unital \(A\), \(K_0(A)\) is the Grothendieck group of \(V(A)\). For nonunital \(A\), set \(K_0(A)={\rm Ker\,} K_0(\epsilon)\). Write \([p]_{K_0}\) for the \(K_0\)-class of a projection, \(K_0(A)^+\) for the image of \(V(A)\), and \(\Sigma(A)=\{[p]_{K_0}:p\text{ a projection in }A\}\) for the scale. In both cases,
\[
 K_1(A)=\directlim_n U(M_n(A^\sim))/U_0(M_n(A^\sim)),
\]
where \(U_0\) denotes the identity component and the maps are induced by \(u\mapsto u\oplus1\). Write \([u]_{K_1}\) for a unitary class.
\end{definition}

\begin{definition}[{\cite[Definition~3.14]{ERR}}]\label{def:swc}
Let \(A\) be a \(C^*\)-algebra. We say that \(A\) has \emph{weak
cancellation} if, for all projections \(p,q\in A\) generating the same
closed ideal \(J\), the equality \([p]_{K_0}=[q]_{K_0}\) in \(K_0(J)\) implies
\(p\sim_{\mathrm{MvN}}q\). We say that \(A\) has \emph{stable weak
cancellation} if \(M_n(A)\) has weak cancellation for every \(n\ge1\).
Equivalently, \(A\otimes\K\) has weak cancellation. A projection in
\(A\) is \emph{full} if the closed ideal it generates is \(A\).
\end{definition}

\begin{definition}[{\cite{Rie,BP}}]
For \(m\ge1\), put
\[
 \operatorname{Lg}_m(B)
 =\{(a_1,\ldots,a_m)\in B^m:
                      \textstyle\sum_{i=1}^m a_i^*a_i\text{ is invertible}\}
\]
for unital \(B\). The stable rank \(\sr(A)\) is the least \(m\ge1\) for which \(\operatorname{Lg}_m(A^\sim)\) is dense in \((A^\sim)^m\). The real rank \({\rm rr}(A)\) is the least \(m\ge0\) for which
\(\operatorname{Lg}_{m+1}(A^\sim)\cap((A^\sim)_{\mathrm{sa}})^{m+1}\)
is dense in \(((A^\sim)_{\mathrm{sa}})^{m+1}\). If no such integer exists, the corresponding rank is \(\infty\).
\end{definition}

\begin{definition}
For \(a,b\in(A\otimes\K)_+\), write \(a\precsim b\) if \(x_j^*bx_j\to a\) in norm for some \(x_j\in A\otimes\K\). Write \([a]\) for the class of \(a\) under mutual subequivalence. The Cuntz semigroup \(\Cu(A)\) consists of these classes, with
\[
 [a]\le[b]\ \Longleftrightarrow\ a\precsim b,\qquad
 [a]+[b]=[a\oplus b].
\]
For projections,
\[
           p\precsim q\quad\Longleftrightarrow\quad
           p\precsim_{\mathrm{MvN}}q .
\]
The map \(V(A)\to\Cu(A)\), \([p]_{\mathrm{MvN}}\mapsto[p]\), preserves addition and comparison. When \(A\) is stably finite, mutual subequivalence of projections implies Murray--von Neumann equivalence, so this map is an order embedding.
\end{definition}

In a Cuntz semigroup, we write \(x\ll y\) if every increasing sequence \((y_n)_{n\ge 1}\) with \(y\le\sup_n y_n\) satisfies \(x\le y_n\) for some \(n\). An element \(x\) is \emph{compact} if \(x\ll x\). The semigroup is \emph{algebraic} if every element is the supremum of an increasing sequence of compact elements, and it has \emph{weak cancellation} if \(x+z\ll y+z\) implies \(x\ll y\) for all elements \(x,y,z\) \cite{RW}.

\begin{definition}[{\cite[Definitions~5.2.16 and~5.6.1]{APT}}]
Let $(S,\leq)$ be a positively ordered monoid. In the following definitions, the implications are required for all $x,y\in S$. We say that $S$ is

(1) \emph{unperforated} if $nx\le ny \text{ for some } n\ge1$ implies $x\le y,$

(2) \emph{nearly unperforated}  if $ nx\le ny \text{ for all sufficiently large } n$ implies $x\le y,$

(3) \emph{almost unperforated} if $ (n+1)x\le ny \text{ for some } n\ge1 $ implies $x\le y .$

\end{definition}

These conditions are listed in decreasing order of strength. Near unperforation is equivalent to \(2x\le2y,\ 3x\le3y\Rightarrow x\le y\) \cite[Lemma~5.6.2 and Proposition~5.6.3]{APT}.

\begin{definition}[{\cite{Bla,LR,NW}}] An \emph{AF algebra} is an inductive limit of finite-dimensional algebras. An \emph{AT algebra} is an inductive limit of finite direct sums of matrix algebras over points and circles. An algebra is \emph{subhomogeneous} if the dimensions of its irreducible representations are uniformly bounded, and \emph{ASH} if it is an inductive limit of subhomogeneous algebras. An algebra \(A\) is \emph{locally subhomogeneous} if, for every finite set \(F\subset A\) and \(\varepsilon>0\), there is a subhomogeneous subalgebra \(S\subset A\) such that \(\dist(a,S)<\varepsilon\) for all \(a\in F\).
\end{definition}

\begin{definition}[{\cite{Busby}}]
An extension of \(B\) by \(I\) is a short exact sequence
\[
 e:\quad 0\longrightarrow I\xrightarrow{\iota}A
                    \xrightarrow{\pi}B\longrightarrow0,
\]
with Busby invariant \(\tau\colon B\to\mathcal M(I)/I\).
The extension is \emph{unital} if \(A\) is unital and \(\tau\) is unital.
It is \emph{norm-full} if \(\tau(b)\) generates \(\mathcal M(I)/I\)
as a closed ideal for every nonzero \(b\in B_+\).
\end{definition}
If \(B\) is unital and simple and \(I\) is nonunital, every unital
extension of \(B\) by \(I\) is essential and norm-full.

\begin{definition}[{\cite{Bla}}]
The six-term exact sequence of $
          0\longrightarrow I\xrightarrow{\iota}A
                    \xrightarrow{\pi}B\longrightarrow0
$ is
\[
\begin{CD}
 K_0(I) @>{K_0(\iota)}>> K_0(A) @>{K_0(\pi)}>> K_0(B)\\
 @A{\delta_1}AA @. @VV{\delta_0}V\\
 K_1(B) @<{K_1(\pi)}<< K_1(A) @<{K_1(\iota)}<< K_1(I) .
\end{CD}
\]
Its boundary maps are the \emph{exponential map} \(\delta_0\) and the \emph{index map} \(\delta_1\). With unitizations and matrix amplifications understood, for a self-adjoint lift \(b\) of a projection \(p\) and a unitary lift \(U\) of \(u\oplus u^*\), the formulas are
\[
 \delta_0([p]_{K_0})=[\exp(2\pi i b)]_{K_1},\qquad
 \delta_1([u]_{K_1})=[UeU^*]_{K_0}-[e]_{K_0},\quad e=\diag(1,0).
\]
If \(u\) has a partial-isometry lift \(v\), then
\[
       \delta_1([u]_{K_1})=[1-v^*v]_{K_0}-[1-vv^*]_{K_0} .
\]
\end{definition}
\begin{remark}\label{rem:swc}
If \(\sr(A)=1\), then \(A\) has cancellation of projections and hence
stable weak cancellation \cite[Proposition~6.5.1]{Bla}. In particular, every AT algebra has stable weak cancellation.
Let \(0\to I\to A\to B\to0\) be an extension with \(\RR(A)=0\).
Then \(A\) has stable weak cancellation if and only if both \(I\) and
\(B\) have stable weak cancellation \cite[Lemma~3.15]{ERR}.
\end{remark}
We use Spielberg's criterion to check stable finiteness from the
index map and the positive cone of the ideal group.
\begin{proposition}[{\cite[Lemma~1.5]{Sp}}]\label{prop:spielberg}
Let \(0\to I\xrightarrow{\iota}A\to B\to0\) be an extension with \(I\) and \(B\) stably finite. Then \(A\) is stably finite if and only if
\[
                 \im\delta_1\cap K_0(I)^+=\{0\}.
\]
\end{proposition}
The boundary maps also determine the permanence of real rank zero
and stable rank one under the corresponding hypotheses on the ideal
and quotient.
\begin{proposition}[{\cite[Proposition~4]{LR}}]\label{prop:lr}
Let
$
             0\rightarrow I\rightarrow A
                         \rightarrow B\rightarrow0
$
be an extension of $C^*$-algebras.
Then the following assertions hold.

{\rm (i)} If \(I,B\) have real rank zero, then \(\RR(A)=0\) if and only if \(\delta_0=0\).

{\rm (ii)} If \(I,B\) have stable rank one, then \(\sr(A)=1\) if and only if \(\delta_1=0\).
\end{proposition}

\section{An Extension Construction}\label{sec:disk}

We construct an extension of a simple AT algebra by a stable simple AF algebra with nonzero index map. The initial disk block contains a partial isometry whose defects give the required failure of cancellation.

\begin{example}\label{ex:firstblock}
Let \(r,s\in\mathbb N\), and let \(\Dsk\) be the closed unit disk. Define
\[
 \alpha_{r,s}\colon M_r\oplus C(\T,M_s)\longrightarrow C(\T,M_{r+s}),
 \qquad
 \alpha_{r,s}(a,h)(\xi)=
      \begin{pmatrix}a&0\\0&h(\xi^s)\end{pmatrix}.
\]
Here \(h(\xi^s)\in M_s\) is obtained by composing \(h\) with the degree-\(s\) circle map. The \(2\)-NCCW complex \(E(r,s)\) is the pullback
\[
\begin{CD}
 E(r,s) @>>> C(\Dsk,M_{r+s})\\
 @VVV @VV{\mathrm{res}_{\T}}V\\
 M_r\oplus C(\T,M_s) @>{\alpha_{r,s}}>> C(\T,M_{r+s}).
\end{CD}
\]
Thus
\[
 E(r,s)=\{(f,a,h):f|_{\T}=\alpha_{r,s}(a,h)\}.
\]

This is a \(2\)-NCCW complex in the sense of \cite{ELP}. The quotient map \(\pi(f,a,h)=h\) is onto by extension of matrix entries from \(\T\) to \(\Dsk\). Its kernel and a splitting of its constant-coordinate map are
\[
 I(r,s)=\{(f,a,0):f|_{\T}=a\oplus0_s\},\qquad
 a\longmapsto(a\oplus0_s,a,0).
\]

Set
\[
 s_n=n!,\qquad r_1=1,\qquad
 r_{n+1}=(n+1)r_n+ns_n,
\]
and take
\[
        E_n=E(r_n,s_n),\qquad I_n=I(r_n,s_n),\qquad
        D_n=C(\T,M_{s_n}).
\]
Write \(\alpha_n=\alpha_{r_n,s_n}\), \(\pi_n(f,a,h)=h\), and
\(\iota_n\colon I_n\hookrightarrow E_n\). These give the extensions
\[
            0\longrightarrow I_n\xrightarrow{\iota_n}E_n
                   \xrightarrow{\pi_n}D_n\longrightarrow0.
\]

Let \((z_j)_{j\ge 1}\) be an enumeration of \((\Q+i\Q)\cap\operatorname{int}\Dsk\). Since every integer \(n\ge 1\) has unique representations
\[
 n=2^m+k \quad (m\ge 0,\ 0\le k<2^m)
 \qquad\text{and}\qquad
 n=\frac{(i+j-2)(i+j-1)}{2}+i \quad (i,j\ge 1),
\]
we may define sequences \((\zeta_n)_{n\ge 1}\) in \(\T\) and \((x_n)_{n\ge 1}\) in \(\operatorname{int}\Dsk\) by
\[
 \zeta_{2^m+k}=e^{2\pi i k/2^m},
 \qquad
 x_n=z_j^{\,s_i/s_n},
\]
where \(z_j^{\,s_i/s_n}\) denotes the principal \((s_n/s_i)\)-th root of \(z_j\) (with \(0^{s_i/s_n}=0\)). Notice that \(n\ge i\) holds for all \(i,j\ge 1\). For each fixed \(i\ge 1\), the indices \(n\) corresponding to \((i,j)\) (\(j\ge 1\)) therefore satisfy \(n\ge i\) and \(x_n^{s_n/s_i}=z_j\), while \(\{\zeta_n : n\ge i\}\) contains all \(2^m\)-th roots of unity whenever \(2^m\ge i\). Consequently,
\begin{equation}\label{eq:points}
 \overline{\{x_n^{s_n/s_i}:n\ge i\}}=\Dsk,
 \qquad
 \overline{\{\zeta_n:n\ge i\}}=\T.
\end{equation}

The connecting maps on the quotient are
\[
 \beta_{n,n+1}\colon D_n\longrightarrow D_{n+1},\qquad
 \beta_{n,n+1}(h)(z)=
       \begin{pmatrix}
          h(z)&0\\0&h(\zeta_n)^{\oplus n}
       \end{pmatrix}.
\]
To define the disk maps, decompose
\(\mathbb C^{r_{n+1}+s_{n+1}}\) initially into blocks of sizes
\[
 r_n,\quad s_n,\quad n(r_n+s_n),\quad ns_n.
\]
Let \(U_n\) interchange the second and third blocks. In these domain and target block decompositions,
\[
 U_n=
 \begin{pmatrix}
  1_{r_n}&0&0&0\\
  0&0&1_{n(r_n+s_n)}&0\\
  0&1_{s_n}&0&0\\
  0&0&0&1_{ns_n}
 \end{pmatrix}.
\]
Define \(\psi_{n,n+1}\colon E_n\to E_{n+1}\) by \(\psi_{n,n+1}(f,a,h)=(\widehat{f} ,\widehat{a},\widehat{h})\), where
\begin{equation*}\label{eq:diskmap}
\begin{gathered}
 \widehat{a}=
 \begin{pmatrix}a&0\\0&f(x_n)^{\oplus n}\end{pmatrix},
 \qquad \widehat{h}=\beta_{n,n+1}(h),\\
 \widehat{f}(\lambda)=
 U_n
 \begin{pmatrix}
  f(\lambda^{n+1})&0&0\\
  0&f(x_n)^{\oplus n}&0\\
  0&0&h(\zeta_n)^{\oplus n}
 \end{pmatrix}
 U_n^*
                   \qquad(\lambda\in\Dsk).
\end{gathered}
\end{equation*}
For \(\xi\in\T\), this gives the boundary identity
\[
 \widehat{f}(\xi)=
 \begin{pmatrix}
  a&0&0&0\\
  0&f(x_n)^{\oplus n}&0&0\\
  0&0&h(\xi^{s_{n+1}})&0\\
  0&0&0&h(\zeta_n)^{\oplus n}
 \end{pmatrix}
 =\begin{pmatrix}\widehat{a}&0\\0&\widehat{h}(\xi^{s_{n+1}})\end{pmatrix}.
\]
Thus \(\psi_{n,n+1}(E_n)\subset E_{n+1}\). Each \(\psi_{n,n+1}\) is unital and injective, since \(\lambda\mapsto\lambda^{n+1}\) is onto and \(\alpha_n\) is injective.

Let \(\lambda_{n,n+1}=\psi_{n,n+1}|_{I_n}\). The diagrams
\[
\begin{CD}
 0 @>>> I_n @>{\iota_n}>> E_n @>{\pi_n}>> D_n @>>>0\\
 @. @VV{\lambda_{n,n+1}}V @VV{\psi_{n,n+1}}V @VV{\beta_{n,n+1}}V @.\\
 0 @>>> I_{n+1} @>{\iota_{n+1}}>> E_{n+1} @>{\pi_{n+1}}>> D_{n+1} @>>>0
\end{CD}
\]
give, on taking limits,
\[
        0\longrightarrow I\xrightarrow{\iota}E
                    \xrightarrow{\pi}D\longrightarrow0,
\]
where \(I=\directlim(I_n,\lambda_{n,n+1})\), \(E=\directlim(E_n,\psi_{n,n+1})\), and
\(D=\directlim(D_n,\beta_{n,n+1})\).
For \(m>n\), we use the notation
\[
 \psi_{n,m}=\psi_{m-1,m}\circ\cdots\circ\psi_{n,n+1},
 \qquad \psi_{n,n}=\mathrm{id}_{E_n},
\]
and define \(\lambda_{n,m}\) and \(\beta_{n,m}\) in the same way. The
canonical maps into the limits are denoted by \(\psi_{n,\infty}\),
\(\lambda_{n,\infty}\) and \(\beta_{n,\infty}\), respectively.
\end{example}
\begin{proposition}\label{prop:ideal}
The ideal \(I\) is stable simple AF, with \(K_0(I)=\Q^2\) and
\[
 K_0(I)^+=\{0\}\cup\{(x,y)\in\Q^2:x+y>0\}.
\]
\end{proposition}
\begin{proof}
Let \(p_1=(\diag(1,0),1,0)\in I_1\), and put
\[
p_n=\lambda_{1,n}(p_1),\qquad p= \lambda_{1,\infty} (p_1).
\]
The disk and constant coordinates of \(p_n\) have rank \(s_n\), so \(p_n\)
is full in \(I_n\). Write \(S^2=\Dsk/\T\), denote quotient points by
\(\dot\lambda\), and put \(\gamma_n(\dot\lambda)=\dot{\lambda^{n+1}}\).
Compression gives
\[
\begin{gathered}
 p_nI_np_n\cong C(S^2,M_{s_n}),\qquad
 B:=pIp\cong\directlim\bigl(C(S^2,M_{s_n}),\rho_{n,n+1}\bigr),\\
 \rho_{n,n+1}(g)=
 \begin{pmatrix}g\circ\gamma_n&0\\0&g(\dot x_n)^{\oplus n}\end{pmatrix}.
\end{gathered}
\]

Use the splitting
\[
 K_0(C(S^2))=\Z[1]_{K_0}\oplus\widetilde K_0(C(S^2))
            =\Z[1]_{K_0}\oplus\Z\eta,
\]
where the sign of \(\eta\) is chosen so that disk restriction sends the
winding-number-one class to \(\eta\) under the index map
\cite[Section~9.3]{Bla}. Through the full corner, write \((a,c)\) for the
class \(a[1]_{K_0}+c\eta\) in \(K_0(I_n)\).
A nonzero projection has positive rank. Conversely, the clutching
construction gives a line bundle on \(S^2\) of each reduced class
\(c\eta\), \(c\in\Z\). Adding \(a-1\) trivial lines represents
\((a,c)\) for every \(a\ge1\). Hence
\[
 K_0(I_n)=\Z^2,
 \qquad K_0(I_n)^+=\{0\}\cup\{(a,c)\in\Z^2:a>0\}.
\]
The variable map has degree \(n+1\), and point evaluations kill the reduced
class. Therefore
\[
 K_0(\lambda_{n,n+1})=
 \begin{pmatrix}n+1&0\\0&n+1\end{pmatrix},
 \qquad
 \begin{pmatrix}a\\c\end{pmatrix}
 \longmapsto\begin{pmatrix}(n+1)a\\(n+1)c\end{pmatrix}.
\]
Consequently,
\[
\begin{CD}
K_0(I_1) @>>> K_0(I_2) @>{}>> K_0(I_3) @>{}>> \cdots @>>>K_0(I)\\
@VV{\cong}V @VV{\cong}V @VV{\cong}V @. @VV{\cong}V@.\\
 \Z^2 @> {\times 2}>> \Z^2 @>{\times 3}>> \Z^2 @>{}>> \cdots  @>>>\Q^2 \,.
\end{CD}
\]
The map \([(a,c),n]\mapsto(a/s_n,c/s_n)\) identifies the ordered limit
with \(\Q^2\) having cone \(\{0\}\cup\{(u,b):u>0\}\): the connecting
maps preserve positive rank, and every rational pair occurs at a common
factorial denominator. We use the coordinates
\begin{equation}\label{eq:ideal-coordinates}
 [(a,c),n]\longmapsto
       \left(\frac{a+c}{s_n},-\frac{c}{s_n}\right).
\end{equation}
In these coordinates the cone is
\(\{0\}\cup\{(x,y):x+y>0\}\), and \([p]_{K_0}=(1,0)\).
We also have
\[
 K_1(I_n)=K_1\bigl(C(S^2,M_{s_n})\bigr)=0,
 \qquad K_1(B)=\directlim K_1(I_n)=0.
\]

By \eqref{eq:points}, every nonzero positive stage element in the system for
\(B\) acquires a nonzero constant block at a later stage. Its image is then
nonzero at every point of \(S^2\), and hence is full. Thus \(B\) is simple.

Since this construction is a slight modification of Goodearl's, an argument similar to that of \cite[Theorem~9]{Goo} shows that \(B\) has real rank zero. (One can also see directly that the connecting maps have arbitrarily small spectral variation, since each composite has only one nonconstant eigenvalue map and the constant evaluation points become dense in \(S^2\).
Then \cite[Corollary~2.5]{EG} gives the same conclusion.)

The ordered group \((K_0(I),K_0(I)^+,(1,0))\) is a simple dimension
group. The AF range theorem gives a simple unital AF algebra \(F_0\)
with this scaled ordered \(K_0\)-group \cite[Theorem~2.2]{EHS}.

Applying Elliott--Gong's classification theorem
\cite[Theorem~5.1]{EG} gives \(B\cong F_0\).

Since \(p\) is full in \(I\), Brown's theorem \cite[Theorem~2.8]{Br} yields \(I\otimes\K\cong B\otimes\K\). As \(B\otimes\K\) is AF and hereditary subalgebras of AF algebras are AF, it follows that \(I\) is AF.

It remains to prove that \(I\) is stable. We show that every positive \(K_0\)-class is represented by a projection in \(I\).  Let \(q_n\in I\) be the image of
\((\diag(1_{r_n},0),1_{r_n},0)\in I_n\). Then
\[
 [q_n]_{K_0}=(\frac{r_n}{s_n},0),\qquad
 \frac{r_{n+1}}{s_{n+1}}-\frac{r_n}{s_n}
                  =\frac n{n+1}\ge\frac12.
\]
In particular, \(r_n/s_n\to\infty\).  Let \(g=(x,y)\ne0\) be positive in \(K_0(I)^+\), so \(x+y>0\), and choose \(n\)
with \(r_n/s_n>x+y\). Then
\[
             [q_n]_{K_0}-g=(r_n/s_n-x,-y)\in K_0(I)^+.
\]
Choose a projection in \(I\otimes\K\) representing \(g\). By comparison of projections in the AF algebra \(I\otimes\K\), this projection is equivalent to a subprojection \(e\) of \(q_n\). Thus \(e\in I\) and \(g=[e]_{K_0}\), so \(\Sigma(I)=K_0(I)^+\).
By Elliott's classification of AF algebras \cite[Theorem~7.3.2]{Bla}, \(I\cong I\otimes\K\).

\end{proof}

\begin{proposition}\label{prop:quotient}
The quotient \(D\) is the Bunce--Deddens algebra of type
\(\prod_{\ell\ \mathrm{prime}}\ell^\infty\), with
\[
 (K_0(D),K_0(D)^+,[1_D]_{K_0})=(\Q,\Q_{\ge0},1),
 \qquad K_1(D)=\Z.
\]

\end{proposition}
\begin{proof}
The inductive system is
\[
C(\T)\xrightarrow{\beta_{1,2}}C(\T,M_2)
          \xrightarrow{\beta_{2,3}}C(\T,M_6)
          \longrightarrow\cdots\longrightarrow D.
\]
The induced maps in \(K\)-theory are
\(K_0(\beta_{n,n+1})=\times(n+1)\) and \(K_1(\beta_{n,n+1})=\mathrm{id}\). Thus
\[
\Z\xrightarrow{\times2}\Z
          \xrightarrow{\times3}\Z\rightarrow\cdots \rightarrow K_0(D)
          \cong\Q,
 \]
 \[
\Z\xrightarrow{\mathrm{id}}\Z
          \xrightarrow{\mathrm{id}}\Z\rightarrow\cdots\rightarrow K_1(D)
          \cong\Z.
\]
Moreover,
\[
 K_0(D)^+=\bigcup_{n\ge1}\{\frac a{n!} \mid a\in \mathbb{Z}, a\geq 0\}=\{x\in\Q \mid x\ge0\},
 \qquad [1_D]_{K_0}=s_n/n!=1.
\]
Since every tail of \((\zeta_n)\) is dense in \(\T\), the algebra \(D\) is simple. Moreover, the connecting maps are of Goodearl's form, and for each fixed \(i\), the weighted identity ratio of \(\beta_{i,m}\) is \(s_i/s_m\), which tends to \(0\) as \(m\to\infty\). Hence \(D\) has real rank zero by \cite[Theorem~9]{Goo}.

By Elliott's classification theorem \cite[Theorem~7.1]{Ell93} and the description in \cite[Example~6.7]{PRRS}, \(D\) is isomorphic to the Bunce--Deddens algebra of type \(\prod_{\ell\ \mathrm{prime}}\ell^\infty\).
\end{proof}

The extension
\[
 0\longrightarrow I_n\xrightarrow{\iota_n}E_n
                     \xrightarrow{\pi_n}D_n\longrightarrow0
\]
induces the six-term exact sequence
\[
\begin{CD}
 \Z^2 @>{K_0(\iota_n)}>> K_0(E_n) @>{K_0(\pi_n)}>> \Z\\
 @A{\left(\begin{smallmatrix}0\\s_n\end{smallmatrix}\right)}AA
                 @. @VV{0}V\\
 \Z @<<< K_1(E_n) @<<< 0.
\end{CD}
\]
The left boundary map sends \(1\) to \((0,s_n)\) in the coordinates
\((a,c)\), since \(\xi\mapsto\xi^{s_n}\) has degree \(s_n\).
Consequently
\[
             K_0(E_n)=\Z^2\oplus\Z/s_n\Z,\qquad K_1(E_n)=0.
\]
Use coordinates \((a,b,\bar c)\), where \(a,b\) are the constant and
quotient ranks, and \(\bar c\in\Z/s_n\Z\) is the image of the second
coordinate of \(K_0(I_n)\). The projection
\((0_{r_n}\oplus e_{11},0,e_{11})\) lifts the quotient rank generator. Then
\[
  K_0(\iota_n)(a,c)=(a,0,\bar c),\qquad
  K_0(\pi_n)(a,b,\bar c)=b,\qquad [1_{E_n}]_{K_0}=(r_n,s_n,\bar0).
\]
The disk map multiplies the second coordinate of \(K_0(I_n)\) by
\(n+1\). The \(n\) disk evaluations contribute \(n(a+b)\) to the
constant rank. Hence
\[
 \Psi_{n,n+1}=K_0(\psi_{n,n+1}),\qquad
 \Psi_{n,n+1}
 \begin{pmatrix}a\\b\\\bar c\end{pmatrix}
 =
 \begin{pmatrix}
      n+1&n&0\\
      0&n+1&0\\
      0&0&n+1
 \end{pmatrix}
 \begin{pmatrix}a\\b\\\bar c\end{pmatrix},
\]
where the last coordinate in the image is \(\overline{(n+1)c}\in\Z/s_{n+1}\Z\).
In particular,
\[
 K_0(E)=\directlim\bigl(\Z^2\oplus\Z/s_n\Z,\Psi_{n,n+1}\bigr),\quad
 K_1(E)=0 .
\]
Naturality gives the commutative diagram with exact rows
\[
\begin{CD}
 \Z @>{\left(\begin{smallmatrix}0\\s_n\end{smallmatrix}\right)}>> \Z^2 @>{K_0(\iota_n)}>>
                  \Z^2\oplus\Z/s_n\Z @>{K_0(\pi_n)}>> \Z @>>>0\\
 @VV{\mathrm{id}}V @VV{(n+1)\mathrm{id}}V
       @VV{\Psi_{n,n+1}}V @VV{\times(n+1)}V @.\\
 \Z @>{\left(\begin{smallmatrix}0\\s_{n+1}\end{smallmatrix}\right)}>> \Z^2 @>{K_0(\iota_{n+1})}>>
                  \Z^2\oplus\Z/s_{n+1}\Z @>{K_0(\pi_{n+1})}>> \Z @>>>0 .
\end{CD}
\]

\begin{proposition}\label{prop:six}
The extension
\[
 0\longrightarrow I\xrightarrow{\iota}E
                  \xrightarrow{\pi}D\longrightarrow0
\]
induces the six-term exact sequence
\begin{equation}\label{eq:six}
\begin{CD}
 \Q^2 @>{K_0(\iota)}>> \Q\oplus(\Q/\Z)\oplus\Q @>{K_0(\pi)}>> \Q\\
 @A{\delta_1}AA @. @VV{0}V\\
 \Z @<<< 0 @<<< 0 .
\end{CD}
\end{equation}
Here
\[
 \delta_1(n)=(n,-n),\qquad
 K_0(\iota)(x,y)=(x+y,x+\Z,0),\qquad
 K_0(\pi)(a,\bar b,c)=c,
\]
and \([1_E]_{K_0}=(0,0,1)\). Moreover, \(E\) has real rank zero, and
\begin{equation*}
\begin{split}
 K_0(E)^+={}&\{0\}\cup
       \{(a,\bar b,0):a\in\Q_{>0},\ \bar b\in\Q/\Z\}\\
       &\quad\cup
       \{(a,\bar b,c):a\in\Q,\ \bar b\in\Q/\Z,\ c\in\Q_{>0}\}.
\end{split}
\end{equation*}
\end{proposition}
\begin{proof}
At stage \(n\), the index class \((0,s_n)\) represents \((0,1)\), and hence corresponds to \((1,-1)\). Passing to the inductive limit yields the exact sequence
\[
 0\longrightarrow\Q^2/\Z(1,-1)\longrightarrow K_0(E)
   \longrightarrow\Q\longrightarrow0,
\]
together with the isomorphism

 \(\Q^2/\Z(1,-1)\cong\Q\oplus(\Q/\Z)\) given by \((x,y)+\Z(1,-1)\mapsto(x+y,x+\Z)\).

 Since the kernel \(\Q\oplus(\Q/\Z)\) is divisible, this extension splits, and there exists a homomorphism \(\Q\to\Q\oplus(\Q/\Z)\) taking any prescribed value at \(1\). Adjusting the splitting accordingly, we may assume that it sends \(1\) to \([1_E]_{K_0}\). In these coordinates, the maps take the stated form.

The ideal $I$ and quotient $D$ have real rank zero by
Propositions~\ref{prop:ideal} and~\ref{prop:quotient}. Since the exponential
map vanishes, Proposition~\ref{prop:lr} gives
\(\RR(E)=0\). The stable simple AF ideal \(I\) has the corona
factorization property \cite[Theorem~3.10\textup{(i)}]{ERR}, and the simple AT algebra \(D\)
has cancellation of projections.  Hence \cite[Proposition~3.13]{ERR} gives the positive cone
\[
 K_0(E)^+=K_0(\iota)(K_0(I)^+)
       \sqcup\{g\in K_0(E):K_0(\pi)(g)>0\}.
\]
Substituting the positive cone of the ideal and the maps above yields the desired result.
\end{proof}
Write \(\Zs\) for the Jiang--Su algebra. An algebra \(A\) is \emph{\(\Zs\)-stable} if \(A\cong A\otimes\Zs\) \cite{JS99,TW}.
\begin{proof}[Proof of Theorem~\ref{thm:A}]
By construction, \(E\) is a unital inductive limit of
\(2\)-NCCW complexes, and Proposition~\ref{prop:six} gives \(\RR(E)=0\). Both \(I\) and \(D\) have stable rank one.
Since \(\delta_1\ne0\), Proposition~\ref{prop:lr} gives \(\sr(E)>1\).

The extension estimate
\cite[Corollary~4.12]{Rie} gives \(\sr(E)=2\).

The infinite-dimensional simple AF algebra \(I\) and the simple AT
quotient \(D\) absorb \(\Zs\) \cite[Corollary~3.1]{TWASH}.
By \cite[Corollary~3.2 and Theorem~4.3]{TW},  \(E\) is \(\Zs\)-stable.

We now exhibit projections in \(E\) to show the failure of cancellation. Let \(z(\xi)=\xi\) for \(\xi\in\T\), and define
\[
 v_0(\lambda)=
 \begin{pmatrix}
  0&\sqrt{1-|\lambda|^2}\\
  0&\lambda
 \end{pmatrix}
 \qquad(\lambda\in\Dsk).
\]
Since \(v_0(\xi)=\diag(0,\xi)\) for all \(\xi\in\T\), the triple \(v_1=(v_0,0,z)\) belongs to \(E_1=E(1,1)\). Moreover, \(v_0(\lambda)^*v_0(\lambda)=\diag(0,1)\) for all \(\lambda\in\Dsk\), so \(v_1\) is a partial isometry with \(\pi_1(v_1)=z\).

Set \(v=\psi_{1,\infty}(v_1)\), and define the projections
\[
 e=1_E-v^*v,\qquad f=1_E-vv^*,\qquad r=v^*v.
\]
Since \(\pi(v)=\beta_{1,\infty}(z)\) is a unitary representing the generator of \(K_1(D)=\Z\), we have \(e,f\in I\). The index-map computation in Proposition~\ref{prop:six} therefore gives
\[
 [e]_{K_0}-[f]_{K_0}
   =\delta_1([\pi(v)]_{K_1})
   =(1,-1)\ne0
   \qquad\text{in }K_0(I).
\]
If \(e\) and \(f\) are Murray--von Neumann equivalent in \(E\), any implementing partial isometry would lie in \(fEe\subseteq I\), implying \([e]_{K_0}=[f]_{K_0}\) in \(K_0(I)\), a contradiction. Hence \(e\not\sim_{\mathrm{MvN}}f\) in \(E\).

On the other hand, we have orthogonal decompositions \(1_E=e+r=f+vv^*\) together with \(r=v^*v\sim_{\mathrm{MvN}}vv^*\). Consequently,
\[
 e\oplus r\sim_{\mathrm{MvN}}1_E
          \sim_{\mathrm{MvN}}f\oplus vv^*
          \sim_{\mathrm{MvN}}f\oplus r,
\]
which shows that cancellation of projections fails in \(E\).

\end{proof}
\begin{remark}\label{rem:direct-sum}
The direct sum \(E\oplus E\) is a unital, \(\Zs\)-stable \(C^*\)-algebra of real rank zero and an inductive limit of \(2\)-NCCW complexes. It has stable rank two and fails cancellation. Set
\[
 t=(0,\overline{1/2},0)\in K_0(E),
 \qquad
 g=(t,[1_E]_{K_0})\in K_0(E\oplus E).
\]
Then \(t\notin K_0(E)^+\) and \(2t=0\). Since the positive cone of \(K_0(E\oplus E)\) is \(K_0(E)^+\times K_0(E)^+\), we obtain
\[
 g\notin K_0(E\oplus E)^+,
 \qquad
 2g=(0,2[1_E]_{K_0})>0.
\]
Hence the implication \(ng>0\Rightarrow g>0\) fails for \(n=2\), so
\(K_0(E\oplus E)\) is not unperforated for the strict order in the
sense of \cite[Theorem~2.4]{PhillipsLim}.

The inductive system \((E_n,\psi_{n,n+1})\), and hence its direct sum
with itself, has slow dimension growth in the sense of Phillips. Thus \cite[Theorem~2.3]{PhillipsLim} yields
\[
 \tau(p)<\tau(q)\quad\text{for all }\tau\in T(E\oplus E)
 \quad\Longrightarrow\quad
 p\precsim_{\mathrm{MvN}}q
\]
for all projections \(p,q\) in matrix algebras over \(E\oplus E\).
\end{remark}
\section{Changing the \(K_0\)-order of the ideal}\label{sec:strict}

We retain the quotient \(D\) from Section~\ref{sec:disk} and change the
order on the \(K_0\)-group of the ideal. Set
\[
 H=\Q^2,\qquad
 H^+=\{0\}\cup\{(x,y)\in H:x>0,\ y>0\},\qquad
 \kappa=(1,-1).
\]
The ideal \(I\) constructed in Section~\ref{sec:disk} has
\[
 K_0(I)=H,\qquad
 K_0(I)^+=\{0\}\cup\{(x,y)\in H:x+y>0\}.
\]
Both cones meet \(\Z\kappa\) only at zero, but their images in
\(H/\Z\kappa\) differ. This difference will give the extension in
Theorem~\ref{thm:B}.

The ordered group \((H,H^+)\) is countable, directed, unperforated and
has the Riesz interpolation property. Every nonzero positive element
is an order unit, so it is a simple dimension group. By the AF range
theorem \cite[Theorem~2.2]{EHS}, choose a simple unital AF algebra \(J\)
with
\[
 (K_0(J),K_0(J)^+,[1_J]_{K_0})=(H,H^+,(3,3)).
\]
Fix a rank-one projection \(e\in\K\). We identify \(K_0(J\otimes\K)\)
with \(H\) through the corner embedding \(a\mapsto a\otimes e\).

\begin{proposition}\label{prop:structure}
Let
\[
 0\longrightarrow J\otimes\K\xrightarrow{\iota}A
     \xrightarrow{\pi}D\longrightarrow0
\]
be a unital extension. Then \(A\) is separable, nuclear,
\(\Zs\)-stable and of real rank zero, and has stable weak cancellation.
Moreover, \(J\otimes\K\) is the unique nonzero proper ideal of \(A\).
\end{proposition}
\begin{proof}
Separability and nuclearity pass from the ideal and quotient to \(A\).
Since \(K_1(J\otimes\K)=0\), the exponential map vanishes, and
Proposition~\ref{prop:lr} gives \(\RR(A)=0\). The algebras \(J\) and
\(D\) absorb \(\Zs\) \cite[Corollary~3.1]{TWASH}. Absorption passes to
stabilization and extensions \cite[Corollary~3.2 and Theorem~4.3]{TW},
so \(A\) is \(\Zs\)-stable. Since the ideal and quotient have stable
rank one, Remark~\ref{rem:swc} gives stable weak cancellation for \(A\).

The Busby map is unital and hence nonzero. Simplicity of \(D\) makes
it injective, so the extension is essential. Since \(J\otimes\K\)
is simple, every nonzero ideal of \(A\) contains \(J\otimes\K\).
Simplicity of \(D\) now gives the last assertion.
\end{proof}

\begin{proposition}\label{prop:coset}
Let
\[
 0\longrightarrow J\otimes\K\xrightarrow{\iota}A
     \xrightarrow{\pi}D\longrightarrow0
\]
be a unital extension, and let \(h\in H\). Then
\[
 K_0(\iota)(h)\in K_0(A)^+
 \quad\Longleftrightarrow\quad
 (h+\im\delta_1)\cap H^+\ne\varnothing.
\]
If \(p,q\in A\otimes\K\) are projections satisfying
\[
 [\pi(p)]_{K_0}=[\pi(q)]_{K_0}\ne0,
 \qquad
 K_0(\iota)(h)=[q]_{K_0}-[p]_{K_0},
\]
where \(\pi\) also denotes its stabilization, then, for every \(n\ge1\),
\[
 p^{\oplus n}\precsim_{\mathrm{MvN}}q^{\oplus n}
 \quad\Longleftrightarrow\quad
 (nh+\im\delta_1)\cap H^+\ne\varnothing.
\]
Moreover, \(p^{\oplus n}\sim_{\mathrm{MvN}}q^{\oplus n}\) if and only if
\(nh\in\im\delta_1\).
\end{proposition}
\begin{proof}
Exactness gives \(\ker K_0(\iota)=\im\delta_1\). Suppose that
\(K_0(\iota)(h)=[d]_{K_0}\) for a projection \(d\in A\otimes\K\).
Then \([\pi(d)]_{K_0}=0\). Since \(D\) is stably finite,
\(\pi(d)=0\), so \(d\in(J\otimes\K)\otimes\K\). Its class in \(H\)
therefore belongs to \((h+\im\delta_1)\cap H^+\). Conversely, a class
in this intersection is represented by a projection in the stabilized
ideal and maps to \(K_0(\iota)(h)\). This proves the first equivalence.

Suppose that \(p^{\oplus n}\precsim_{\mathrm{MvN}}q^{\oplus n}\).
Choose a partial isometry \(v\) with \(v^*v=p^{\oplus n}\) and
\(vv^*\le q^{\oplus n}\). The projection \(d=q^{\oplus n}-vv^*\)
satisfies
\[
 [d]_{K_0}
 =n\bigl([q]_{K_0}-[p]_{K_0}\bigr)=K_0(\iota)(nh).
\]
The first equivalence gives \((nh+\im\delta_1)\cap H^+\ne\varnothing\).

Conversely, represent a class in \((nh+\im\delta_1)\cap H^+\) by a
projection \(d\) in the stabilized ideal. Then
\[
 [p^{\oplus n}\oplus d]_{K_0}=[q^{\oplus n}]_{K_0}.
\]
By Proposition~\ref{prop:structure}, every projection with nonzero
quotient image is full. Thus these two projections are full, and stable
weak cancellation gives
\[
 p^{\oplus n}\oplus d\sim_{\mathrm{MvN}}q^{\oplus n}.
\]
In particular, \(p^{\oplus n}\precsim_{\mathrm{MvN}}q^{\oplus n}\).
Finally, \(nh\in\ker K_0(\iota)\) is equivalent to equality of the
\(K_0(A)\)-classes of \(p^{\oplus n}\) and \(q^{\oplus n}\).
Stable weak cancellation for these full projections proves the last
assertion.
\end{proof}

\begin{proposition}\label{prop:strictsix}
There exists a unital \(C^*\)-algebra \(\widetilde E\) fitting into an
essential unital extension
\[
 0\longrightarrow J\otimes\K\xrightarrow{\iota}\widetilde E
     \xrightarrow{\pi}D\longrightarrow0.
\]
Under the group identifications in \eqref{eq:six}, all six maps agree
with those for \(E\), and
\[
 [1_{\widetilde E}]_{K_0}=[1_E]_{K_0}=(0,0,1).
\]
The positive cones of the middle algebras differ; more precisely,
\[
 (1,0,0)\in K_0(E)^+\setminus K_0(\widetilde E)^+.
\]
\end{proposition}
\begin{proof}
The quotient \(D\) is separable, nuclear, unital and satisfies the
universal coefficient theorem (UCT), and \(J\otimes\K\) is separable
and stable. Divisibility of \(H\) gives
\[
 \Ext^1_{\Z}(\Q,H)=0,\qquad
 \{f(1):f\in\Hom(\Q,H)\}=H.
\]
Since \(K_1(J\otimes\K)=0\),
\cite[Remark~4.8 and Theorem~4.14]{GR} identify the strong unital
extension group with
\[
 \Ext_{\mathrm{us}}^{-1}(D,J\otimes\K)
 \cong\Hom(K_1(D),K_0(J\otimes\K))=\Hom(\Z,H).
\]
Choose a unital extension whose index map is
\(\delta_1(n)=n\kappa\), and denote its middle algebra by
\(\widetilde E\). The extension is essential by
Proposition~\ref{prop:structure}.

The exponential map vanishes. The calculation in
Proposition~\ref{prop:six} uses only the underlying groups and the
index map, so it gives the same six-term sequence and the same unit
class for \(\widetilde E\).

For \(E\), the class \((0,1)\) is positive in the ideal group and maps
to \((1,0,0)\). For \(\widetilde E\), however, the coset
\[
 (0,1)+\im\delta_1=\{(j,1-j):j\in\Z\}
\]
contains neither zero nor a pair with both coordinates positive.
Proposition~\ref{prop:coset} therefore gives the asserted difference
between the positive cones.
\end{proof}

\begin{proof}[Proof of Theorem~\ref{thm:B}]
Let \(\widetilde E\) be the algebra from
Proposition~\ref{prop:strictsix}. Proposition~\ref{prop:structure}
gives separability, nuclearity, \(\Zs\)-stability and real rank zero.
Since \(\Z\kappa\cap H^+=\{0\}\),
Proposition~\ref{prop:spielberg} gives stable finiteness.
In the corner \(J\otimes e\subset J\otimes\K\), choose projections
\(e_0,f_0\) with
\[
 [e_0]_{K_0}=(1,2),\qquad [f_0]_{K_0}=(1,1).
\]
Set \(p=1_{\widetilde E}-e_0\) and \(q=1_{\widetilde E}-f_0\).
Both projections lift \(1_D\), and
\(K_0(\iota)(h)=[q]_{K_0}-[p]_{K_0}\) for \(h=(0,1)\).
For every \(n\ge1\),
\[
 nh+\im\delta_1=\{(j,n-j):j\in\Z\}.
\]
When \(n=1\), no integer \(j\) satisfies \(0<j<1\), so this coset
does not meet \(H^+\). When \(n\ge2\), it contains the positive
class \((1,n-1)\). Proposition~\ref{prop:coset} gives
\[
 p\not\precsim_{\mathrm{MvN}}q,
 \qquad p^{\oplus n}\precsim_{\mathrm{MvN}}q^{\oplus n}
 \quad(n\ge2).
\]
Thus neither \(V(\widetilde E)\) nor \(\Cu(\widetilde E)\) is nearly
unperforated.
\end{proof}

\begin{corollary}
The semigroup \(\Cu(\widetilde E)\) is algebraic and almost unperforated, and it fails weak cancellation.
\end{corollary}
\begin{proof}
Since \(\widetilde E\) has real rank zero,
\(\Cu(\widetilde E)\) is algebraic
\cite[Remark~5.5.2(2)]{APT}.
Since \(\widetilde E\) is also \(\mathcal Z\)-stable,
\(\Cu(\widetilde E)\) is almost unperforated
\cite[Theorem~4.5]{Ror04}.
Furthermore, \(\Cu(\widetilde E)\) satisfies almost algebraic
order, that is, axiom~(O5), by \cite[Proposition~4.6]{APT}.

If \(\Cu(\widetilde E)\) had weak cancellation, these
properties and \cite[Proposition~5.6.12]{APT} would imply
that it is nearly unperforated, contradicting
Theorem~\ref{thm:B}.
Thus \(\Cu(\widetilde E)\) does not have weak cancellation.
\end{proof}

\begin{remark}\label{rem:amplifications}
Fix an integer \(m\ge2\), and retain \(f_0\) and
\(q=1_{\widetilde E}-f_0\) from the proof of
Theorem~\ref{thm:B}. Choose a projection \(e_m\in J\otimes e\)
with
\[
 [e_m]_{K_0}=\left(1,1+\frac{1}{m-1}\right),
\]
and put \(p_m=1_{\widetilde E}-e_m\).
For \(n\ge1\), the coset appearing in
Proposition~\ref{prop:coset} is
\[
 n\bigl([e_m]_{K_0}-[f_0]_{K_0}\bigr)
       +\im\delta_1
 =
 \left\{\left(j,\frac{n}{m-1}-j\right):j\in\mathbb Z\right\}.
\]
This coset meets \(H^+\) precisely when
\(0<j<n/(m-1)\) for some integer \(j\), which is equivalent
to \(n\ge m\). Hence
\[
 p_m^{\oplus n}\precsim_{\mathrm{MvN}}q^{\oplus n}
 \quad\Longleftrightarrow\quad n\ge m.
\]

\end{remark}

\section{Strong quasidiagonality and local subhomogeneity}\label{sec:local}

\begin{definition}[{\cite{Had,BrownQD}}]
A representation \(\rho\colon A\to\mathcal B(\mathcal H)\) on a separable Hilbert space is \emph{quasidiagonal} if there is an increasing sequence \((P_n)_{n\ge1}\) of finite-rank projections on \(\mathcal H\) which converges strongly to \(1_{\mathcal H}\) and satisfies
\[
                    \|[P_n,\rho(a)]\|\longrightarrow0\qquad(a\in A).
\]
A separable \(C^*\)-algebra is \emph{quasidiagonal} if it admits a faithful quasidiagonal representation. It is \emph{strongly quasidiagonal} (SQD) if every representation on a separable Hilbert space is quasidiagonal.
\end{definition}

A representation is \emph{essential} when its image contains no nonzero compact operators. Voiculescu's theorem gives quasidiagonality of every faithful essential representation of a separable quasidiagonal algebra; see \cite[Definition~2.8 and Theorem~3.11]{BrownQD}.

We record the following standard consequence of Hadwin's finite-set criterion and his result for CCR algebras \cite[Lemma~1 and Proposition~8(2)]{Had}.
\begin{proposition}\label{prop:local-sqd}
Every separable locally subhomogeneous \(C^*\)-algebra is strongly quasidiagonal. In particular, every separable ASH algebra is strongly quasidiagonal.
\end{proposition}
The following consequence of Brown--Dadarlat's extension theorem will be useful. Eckhardt established strong quasidiagonality for certain
extensions with quotient \(C(\mathbb T)\)
\cite[Lemma~2.10 and Theorem~2.12]{Eck}.
We consider here extensions with locally subhomogeneous
quotients.
\begin{proposition}\label{prop:sqd}
Let
\[
0\longrightarrow I\xrightarrow{\iota}A
 \xrightarrow{\pi}B\longrightarrow0
\]
be an essential extension of separable \(C^*\)-algebras.
Suppose that \(I\) is stable, simple and AF, with
\(I\not\cong\mathcal K\), and that \(B\) is locally subhomogeneous.
Then \(A\) is strongly quasidiagonal if and only if it is stably finite.
\end{proposition}

\begin{proof}
The forward implication is immediate, since strong quasidiagonality implies quasidiagonality and every quasidiagonal \(C^*\)-algebra is stably finite \cite[Proposition~3.19]{BrownQD}.

Conversely, suppose that \(A\) is stably finite.
The quotient \(B\) is nuclear and satisfies the UCT
\cite[Theorem~1.1]{DadUCT}; it is strongly quasidiagonal by
Proposition~\ref{prop:local-sqd}.
The AF ideal \(I\) has the \(K_0\)-embedding property
\cite[p.~76]{NBD}. Since \(I\) is stable, it follows from
\cite[Definition~4.3 and Proposition~4.6]{NBD} that \(A\) is
quasidiagonal.

Essentiality and simplicity of \(I\) imply that every nonzero ideal
of \(A\) contains \(I\). It suffices to consider nondegenerate
representations.
Let \(\rho\colon A\to\mathcal B(\mathcal H)\) be a nondegenerate
representation on a separable Hilbert space.
If \(\rho(I)=0\), then \(\rho\) factors through the quotient
\(B\). Since \(B\) is strongly quasidiagonal, \(\rho\) is
quasidiagonal.

Suppose that \(\rho(I)\ne0\). Since every nonzero ideal of
\(A\) contains \(I\), we have \(\ker\rho=0\).
If \(\rho(A)\cap\mathcal K(\mathcal H)\ne\{0\}\), then the
nonzero ideal \(\rho^{-1}(\mathcal K(\mathcal H))\) contains \(I\).
Choose a nonzero projection \(p\in I\). Then \(\rho(p)\) is a compact projection and therefore has
finite rank. Faithfulness of \(\rho\) gives an embedding
\[
 pIp\longrightarrow
 \rho(p)\mathcal K(\mathcal H)\rho(p).
\]
Hence \(pIp\) is finite-dimensional. Since \(I\) is simple,
\(pIp\) is a full matrix algebra. Fullness of \(p\) and
stability of \(I\) now give
\[
 I\cong pIp\otimes\mathcal K\cong\mathcal K,
\]
which is a contradiction. Therefore \(\rho(A)\cap\mathcal K(\mathcal H)=\{0\}\).
Since \(A\) is quasidiagonal, Voiculescu's theorem implies that
\(\rho\) is quasidiagonal \cite[Theorem~3.11]{BrownQD}.
Thus \(A\) is
strongly quasidiagonal.
\end{proof}

\begin{lemma}\label{lem:bidual}
Let \(B\) be a unital locally subhomogeneous \(C^*\)-algebra,
and let \(p,q\in M_n(B)\) be projections.
If \([p]_{K_0}=[q]_{K_0}\) in \(K_0(B)\), then
\(p\sim_{\mathrm{MvN}}q\) in \(M_n(B^{**})\).
\end{lemma}
\begin{proof}
Choose a projection \(r\in M_k(B)\) and
\(v\in M_{n+k}(B)\) such that
\[
 v^*v=p\oplus r,\qquad vv^*=q\oplus r.
\]
By local subhomogeneity, for any sufficiently small
\(\varepsilon>0\), there are a subhomogeneous subalgebra
\(S\subseteq B\), self-adjoint elements
\(a,b\in M_n(S)\), \(c\in M_k(S)\), and
\(y\in M_{n+k}(S)\) such that
\[
 \|a-p\|<\varepsilon,\quad
 \|b-q\|<\varepsilon,\quad
 \|c-r\|<\varepsilon,\quad
 \|y-v\|<\varepsilon.
\]
For \(0<\varepsilon<1/4\), the spectra of \(a,b,c\) avoid
\(1/2\). Continuous functional calculus therefore gives
projections
\[
 p'=\chi_{(1/2,\infty)}(a),\qquad
 q'=\chi_{(1/2,\infty)}(b),\qquad
 r'=\chi_{(1/2,\infty)}(c)
\]
over \(S\), with
\[
 \|p'-p\|<2\varepsilon,\qquad
 \|q'-q\|<2\varepsilon,\qquad
 \|r'-r\|<2\varepsilon.
\]
Set \(x=(q'\oplus r')y(p'\oplus r')\).
For sufficiently small \(\varepsilon\),
\[
 \|x^*x-(p'\oplus r')\|<1,\qquad
 \|xx^*-(q'\oplus r')\|<1.
\]
Thus \(x^*x\) and \(xx^*\) are invertible in their respective
corners. The partial isometry
\(w=x(x^*x)^{-1/2}\) consequently satisfies
\[
 w^*w=p'\oplus r',\qquad ww^*=q'\oplus r'.
\]

Since \(S^{**}\) is finite, cancellation gives
\(p'\sim_{\mathrm{MvN}}q'\) in \(M_n(S^{**})\).
The inclusion \(S\subseteq B\) extends to an injective
\(*\)-homomorphism \(S^{**}\to B^{**}\).
Thus \(p'\sim_{\mathrm{MvN}}q'\) also holds in \(M_n(B^{**})\).

Since \(\|p-p'\|<1\) and \(\|q-q'\|<1\), the standard
result for close projections gives
\[
 p\sim_{\mathrm{MvN}}p',
 \qquad q'\sim_{\mathrm{MvN}}q
 \qquad\text{in }M_n(B).
\]
These equivalences also hold in \(M_n(B^{**})\).
By transitivity,
\[
 p\sim_{\mathrm{MvN}}p'
  \sim_{\mathrm{MvN}}q'
  \sim_{\mathrm{MvN}}q
 \qquad\text{in }M_n(B^{**}).
\]
\end{proof}
Let $F$ be a unital $C^*$-algebra with  $T(F)\neq \varnothing$. For $\tau\in T(F)$, denote by
\[
 (\tau\otimes\operatorname{Tr})_*
 \colon K_0(F\otimes\mathcal K)\longrightarrow\mathbb R
\]
the homomorphism determined by
\[
 (\tau\otimes\operatorname{Tr})_*
 \bigl([p]_{K_0}-[q]_{K_0}\bigr)
 =
 (\tau\otimes\operatorname{Tr})(p)
 -
 (\tau\otimes\operatorname{Tr})(q),
\]
where $p,q$ are projections in a common finite matrix
corner $M_m(F)\subseteq F\otimes\mathcal K$  and  \(\Tr\) is the unnormalized trace on \(\K\).

\begin{corollary}\label{cor:trace-index}
Let \(A\) be a unital locally subhomogeneous \(C^*\)-algebra,
and let \(I\subseteq A\) be an ideal isomorphic to
\(F\otimes\K\), where \(F\) is a unital \(C^*\)-algebra. Then
\[
 (\tau\otimes\Tr)_*(\delta_1 (x))=0,
       \qquad x\in K_1(A/I),\ \tau\in T(F).
\]
\end{corollary}
\begin{proof}
Since $F$ is isomorphic to a unital corner of $A$, Proposition~\ref{prop:local-sqd} and \cite[Proposition~7.1]{BrownQD} imply that $F$ is quasidiagonal, whence $T(F)\ne\varnothing$ by \cite[Proposition~6.1]{BrownQD}.

Fix \(x\in K_1(A/I)\) and \(\tau\in T(F)\), and put
\(\sigma=\tau\otimes\Tr\).
Using \(K_0(I)\cong K_0(F)\), choose projections
\(p,q\in M_m(F)\subseteq I\) such that
\[
 \delta_1(x)=[p]_{K_0(I)}-[q]_{K_0(I)}.
\]
The map on $K_0$ induced by the embedding $I\hookrightarrow A$
sends this class to zero, so $[p]_{K_0(A)}=[q]_{K_0(A)}$.
By Lemma~\ref{lem:bidual}, there is a partial isometry
\(v\in A^{**}\) satisfying
\[
 v^*v=p,\qquad vv^*=q.
\]
We identify \(I^{**}\) with its canonical weak*-closed
ideal in \(A^{**}\). Since \(p,q\in I^{**}\) and
\(v=qvp\), it follows that \(v\in I^{**}\).

By \cite[Proposition~5.2]{ORT}, \(\sigma\) extends to a
normal tracial weight \(\widetilde\sigma\) on \(I^{**}\).
Therefore
\[
 \sigma(p)
 =\widetilde\sigma(v^*v)
 =\widetilde\sigma(vv^*)
 =\sigma(q).
\]
Both values are finite, so
\[
 (\tau\otimes\Tr)_*(\delta_1(x))
 =\sigma(p)-\sigma(q)=0.
\]
\end{proof}

\begin{proof}[Proof of Theorem~\ref{thm:C}]
Let \(\widetilde E\) be the algebra from
Proposition~\ref{prop:strictsix}. Theorem~\ref{thm:B} and
Proposition~\ref{prop:sqd} show that it is separable, unital, nuclear,
stably finite, \(\Zs\)-stable, of real rank zero and strongly
quasidiagonal.

Define \(s\colon H\to\mathbb R\) by \(s(x,y)=x/3\).
Since \(s(H^+)\subseteq[0,\infty)\) and \(s_1(3,3)=1\),
it is a state on
\((H,H^+,(3,3))\).
As \(J\) is a unital AF algebra, the trace--state correspondence
\cite[p.~51]{Bla} gives a tracial state
\(\tau_1\) on \(J\) satisfying
\[
 (\tau_1\otimes\operatorname{Tr})(p)
 =s_1([p]_{K_0})
 \qquad\text{for every projection }p\in M_m(J).
\]

Since \(\delta_1(1)=(1,-1)\), we have
\[
 (\tau_1\otimes\Tr)_*(\delta_1(1))=\frac13\ne0.
\]
Corollary~\ref{cor:trace-index} therefore shows that
\(\widetilde E\) is not locally subhomogeneous.
\end{proof}

\begin{proposition}\label{prop Irs}
Let \(A\) be a \(C^*\)-algebra, and let \(p,q\in A\otimes\K\) be
projections. If
\(p\otimes1_{\Zs}\precsim_{\mathrm{MvN}}q\otimes1_{\Zs}\), then there exists an integer $k\geq 1$ such that 
$
 n[p]\leq n[q]$ for any $n\geq k$.
\end{proposition}
\begin{proof}
For relatively prime integers \(r,s\ge2\), write
\[
 I_{r,s}=\{f\in C([0,1],M_r\otimes M_s):
 f(0)\in M_r\otimes1_s,\ f(1)\in1_r\otimes M_s\}.
\]
Write \(\Zs=\varinjlim I_{r_n,s_n}\), where the connecting maps are unital and injective.  Since \(p\otimes1_{\Zs}\precsim_{\mathrm{MvN}}q\otimes1_{\Zs}\),
we can approximate them at a sufficiently large stage, then we obtain
\[
 p\otimes1_{I_{r,s}}
 \precsim_{\mathrm{MvN}}
 q\otimes1_{I_{r,s}}
\]
for some relatively prime integers \(r,s\).
Applying the endpoint evaluations
\(I_{r,s}\to M_r\) and \(I_{r,s}\to M_s\) gives
\[
 r[p]\le r[q],\qquad s[p]\le s[q].
\]
Since \(r\) and \(s\) are relatively prime, every sufficiently large
integer is of the form \(ar+bs\), with \(a,b\ge0\).
Additivity therefore gives \(n[p]\le n[q]\) for all sufficiently
large \(n\).
\end{proof}

The following proposition gives the comparison needed for
Theorem~\ref{thm:D}.

\begin{proposition}\label{prop:Ztest}
Let \(A\) be a locally subhomogeneous \(C^*\)-algebra, and let
\(a,b\in(A\otimes\K)_+\). Suppose that
$ 2[a]\le 2[b]$ and $3[a]\le3 [b]$.
Then \(a\otimes1_{\Zs}\precsim b\otimes1_{\Zs}\).
\end{proposition}
\begin{proof}
We first consider projections \(a=p\) and
\(b=q\), with \(p\ne0\). After placing \(p,q\) in orthogonal matrix corners, apply the
perturbation argument in the proof of Lemma~\ref{lem:bidual}
simultaneously to the two comparisons. We obtain orthogonal
projections \(p',q'\) in a finitely generated subhomogeneous
subalgebra \(R\subseteq A\otimes\K\) such that
\[
 1_R=p'+q',\qquad
 p\sim_{\mathrm{MvN}}p',\qquad q\sim_{\mathrm{MvN}}q'.
\]
and
\[
 p'^{\oplus j}\precsim_{\mathrm{MvN}}q'^{\oplus j}
 \quad\text{in }M_j(R),\qquad j=2,3.
\]
Here \(R\) is generated by \(p',q'\) and the matrix entries
implementing these comparisons, all taken in the corner with
unit \(p'+q'\).

Choose a projection \(e\le q'^{\oplus 2}\) such that
\(
 q'^{\oplus 2}-e\sim_{\mathrm{MvN}}p'^{\oplus 2},
\)
let \(L\) be the ideal of \(R\) generated by the entries of \(e\). Denote $Q=R/L$
and let \(\pi\colon R\to  Q\) be the quotient map.

For a representation \(\rho\colon R\to M_m\), denote its matrix
amplification by
\(
 \rho^{(k)}\colon M_k(R)\to M_{km}.
\)
Then
\[
 \rank\rho^{(k)}(e)
 =k\bigl(\rank\rho(q')-\rank\rho(p')\bigr).
\]
As \(L\) is generated by the matrix entries of \(e\),
\[
 \rho(L)=0
 \quad\Longleftrightarrow\quad
 \rho^{(k)}(e)=0
 \quad\Longleftrightarrow\quad
 \rank\rho(p')=\rank\rho(q').
\]

Write \(\bar x=\pi(x)\) for \(x\in R\), and use bars entrywise
for matrices over \(R\).
For \(j=2,3\), choose a projection \(e_j\le q'^{\oplus j}\)
such that
\[
 q'^{\oplus j}-e_j\sim_{\mathrm{MvN}}p'^{\oplus j},
 \qquad e_2=e.
\]
Let \(\sigma\) be an arbitrary irreducible
representation of \(Q\), then
\(\rho=\sigma\circ\pi\) annihilates \(L\), so
\[
 \rank\sigma(\bar p')=\rank\sigma(\bar q'),
 \qquad
 \rank\sigma^{(j)}(\bar e_j)=0.
\]
Thus every matrix entry of \(\bar e_j\) is annihilated by $\sigma$. These representations separate points, so
\(\bar e_j=0\). Consequently,
\[
 \bar p'^{\oplus j}\sim_{\mathrm{MvN}}\bar q'^{\oplus j}
 \qquad(j=2,3).
\]
Letting \(h=[\bar p']_{K_0(Q)}-[\bar q']_{K_0(Q)}\), we obtain
\[
 2h=3h=0,\qquad h=3h-2h=0.
\]

By \cite[Theorem~1.5]{NW} and
\cite[Theorem~2.16]{Phillips}, \(R\) admits a recursive
subhomogeneous decomposition with compact base spaces
\(X_1,\ldots,X_\ell\) of finite covering dimension. Set
\[
 d=\max_{1\le i\le\ell}\dim X_i.
\]  Choose relatively prime integers \(r,s>d+1\), and recall \(I_{r,s}\)  in Proposition~\ref{prop Irs}.
If \(Q\ne0\), then \(Q\otimes I_{r,s}\) has a recursive
subhomogeneous decomposition of dimension at most \(d+1\); see
\cite[Propositions~3.1 and~3.4]{Phillips}.

For every irreducible representation $\sigma: Q\to M_m$, the equality of ranks and \(\bar p'+\bar q'=1_{Q}\) give
\[
 \rank\sigma(\bar p')
 =\rank\sigma(\bar q')
 =\frac{m}{2}\ge1.
\]
Hence, at every point of the product decomposition of
\(Q\otimes I_{r,s}\), both \(\bar p'\otimes1\) and
\(\bar q'\otimes1\) have rank at least
\[
 \min\{r,s\}>\frac{d+1}{2}.
\]
Their \(K_0\)-classes agree, so cancellation
\cite[Proposition~4.3(2)]{Phillips} gives
\[
 \bar p'\otimes1_{I_{r,s}}
 \sim_{\mathrm{MvN}}
 \bar q'\otimes1_{I_{r,s}}.
\]

On each base space of \(R\otimes I_{r,s}\), the set of points
where \(p'\otimes1\) and \(q'\otimes1\) have equal rank is
closed and open. At these points, the evaluation maps factor
through \(Q\otimes I_{r,s}\), so we use the fixed quotient
equivalence.  On the remaining closed and open part, the rank difference is
at least \(\min\{r,s\}>d+1\). Thus
\cite[Proposition~4.2(1)]{Phillips} extends the comparison
from its intersection with the attaching subset.
Proceeding through the pullback decomposition, we preserve
the fixed quotient equivalence on the equal-rank parts
and the values already defined on the attaching subsets.
We obtain
\[
 p'\otimes1_{I_{r,s}}
 \precsim_{\mathrm{MvN}}
 q'\otimes1_{I_{r,s}}.
\]
Together with \(p\sim_{\mathrm{MvN}}p'\) and
\(q\sim_{\mathrm{MvN}}q'\), this proves the same comparison
for \(p,q\). Notice that it holds for all sufficiently large
relatively prime \(r,s\).
The unital embedding \(I_{r,s}\hookrightarrow\Zs\)
(\cite[Proposition~3.3]{RW}) gives
\[p\otimes1_{\Zs}\precsim_{\mathrm{MvN}} q\otimes1_{\Zs}.\]

Now let \(a,b\in (A\otimes\K)_+\) and fix \(\varepsilon>0\). For \(j=2,3\), choose \(x_j\in M_j(A\otimes\K)\) such that
\[
 \|x_j^*b^{\oplus j}x_j-a^{\oplus j}\|<\varepsilon/16.
\]
First, choose \(\delta>0\) small enough that \[\|x_j^*((b-2\delta)_+)^{\oplus j}x_j-a^{\oplus j}\|<\varepsilon/8,\qquad j=2,3.\]
 Next, by local subhomogeneity and continuity of the functional calculus, we obtain a subhomogeneous \(C^*\)-subalgebra \(\widetilde{S}\subseteq A\otimes\K\), positive elements \(a_0,b_0\in \widetilde{S}\), and elements \(y_j\in M_j(\widetilde{S})\) (\(j=2,3\)) satisfying
\[
 \|a-a_0\|<\varepsilon/8,\qquad
 \|b-b_0\|<\delta,\qquad
 \|y_j^*(b_0-2\delta)_+^{\oplus j}y_j-a_0^{\oplus j}\|
 <\varepsilon/4.
\]
Set
\[
 a_1=(a_0-\varepsilon/4)_+,\qquad
 a_2=(a_1-\varepsilon/4)_+,\qquad
 b_1=(b_0-\delta)_+.
\]
Let \(\widetilde{R}\) be generated by \(a_0,b_0\) and the entries of
\(y_2,y_3\), adjoining a unit if necessary.
The cutdown estimates give
\[
 (a-\varepsilon)_+\precsim a_2,\qquad b_1\precsim b,\qquad
 j[a_1]\le j[(b_1-\delta)_+]\quad(j=2,3),
\]
where the last comparisons hold in \(\Cu(\widetilde{R})\).
Again, \(\widetilde{R}\) has a recursive subhomogeneous decomposition of
some finite dimension \(d\).

For \(t>0\), define \(f_t(\lambda)=\lambda(t-\lambda)_+\)
for \(\lambda\ge0\). Let \(\widetilde{L}\) be the ideal of \(\widetilde{R}\) generated by
\(f_{\varepsilon/4}(a_1)\) and \(f_\delta(b_1)\), and let
\(\pi\colon \widetilde{R}\to \widetilde{Q}=\widetilde{R}/\widetilde{L}\) be the quotient map.
Write \(\bar x=\pi(x)\) for \(x\in R\).
Since
\[
 f_{\varepsilon/4}(\bar a_1)=0,\qquad
 f_\delta(\bar b_1)=0.
\]
It follows that
\[
 \operatorname{sp}(\bar a_1)
 \subseteq \{0\}\cup[\varepsilon/4,\infty),
 \qquad
 \operatorname{sp}(\bar b_1)
 \subseteq \{0\}\cup[\delta,\infty).
\]
Hence
\[
 p_0=\chi_{(0,\infty)}(\bar a_1),\qquad
 q_0=\chi_{(0,\infty)}(\bar b_1)
\]
are projections in \(Q\).
The element \(\bar a_1\) is Cuntz equivalent to \(p_0\), and
\((\bar b_1-\delta)_+\precsim q_0\).
Applying the quotient map to the two comparisons above therefore gives
\[
 p_0^{\oplus j}\precsim_{\mathrm{MvN}}q_0^{\oplus j},
 \qquad j=2,3.
\]

The projection argument works for all sufficiently large
relatively prime endpoint sizes. We may therefore choose
relatively prime integers \(r,s>d+1\) such that
\[
 p_0\otimes1_{I_{r,s}}
 \precsim_{\mathrm{MvN}}
 q_0\otimes1_{I_{r,s}}.
\]
Since \(\bar a_2\in p_0\widetilde{Q}p_0\) and
\(\operatorname{Her}(\bar b_1)=q_0\widetilde{Q}q_0\), this gives an element
\(z_0\in \widetilde{Q}\otimes I_{r,s}\) satisfying
\[
 z_0^*z_0=\bar a_2\otimes1,\qquad
 z_0z_0^*\in\operatorname{Her}(\bar b_1\otimes1),
\]
where \(\operatorname{Her}(c)\) denotes the hereditary
subalgebra generated by \(c\).
If \(\widetilde{Q}=0\), take \(z_0=0\) and any relatively prime \(r,s>d+1\).

Let \(\rho\colon \widetilde{R}\to M_m\) be any finite-dimensional
representation. 
The comparison \(2[a_1]\le 2[(b_1-\delta)_+]\) gives
\[
 2\,\rank\rho(a_1)
 \le 2\,\rank\rho((b_1-\delta)_+).
\]
Dividing by \(2\) and using the definitions of the cutdowns,
we obtain
\[
 \rank\rho(a_2)
 \le\rank\rho(a_1)
 \le\rank\rho((b_1-\delta)_+)
 \le\rank\rho(b_1).
\]

Since \(\widetilde{L}\) is generated by
\(f_{\varepsilon/4}(a_1)\) and \(f_\delta(b_1)\),
the condition \(\rho(\widetilde{L})\ne0\) implies that
at least one of their images under \(\rho\) is nonzero.

If \(\rho(f_{\varepsilon/4}(a_1))\ne0\), then \(\rho(a_1)\)
has an eigenvalue in \((0,\varepsilon/4)\). This eigenvalue
becomes zero in
\(\rho(a_2)=(\rho(a_1)-\varepsilon/4)_+\), so
\[
 \rank\rho(a_2)<\rank\rho(a_1).
\]
Similarly, if \(\rho(f_\delta(b_1))\ne0\), then \(\rho(b_1)\)
has an eigenvalue in \((0,\delta)\), and therefore
\[
 \rank\rho((b_1-\delta)_+)<\rank\rho(b_1).
\]
Thus at least one inequality in the preceding chain is strict.
Since matrix ranks are integers, we conclude that
\[
 \rank\rho(b_1)\ge\rank\rho(a_2)+1
 \qquad\text{whenever }\rho(\widetilde{L})\ne0.
\]

Use the product recursive subhomogeneous decomposition of
\(\widetilde{R}\otimes I_{r,s}\), whose base spaces have dimension at most
\(d+1\)(\cite[Proposition~3.4]{Phillips}).
Since \(I_{r,s}\) is nuclear, we identify
\[
 (\widetilde{R}\otimes I_{r,s})/(\widetilde{L}\otimes I_{r,s})
 \cong \widetilde{Q}\otimes I_{r,s}.
\]
On each base space \(X\), let \(Y\) be the set of points at
which every element of \(\widetilde{L}\otimes I_{r,s}\) evaluates to zero.
This set is closed by continuity.
At points of \(Y\), evaluation factors through
\(\widetilde{Q}\otimes I_{r,s}\), so \(z_0\) prescribes a continuous
section there.

At every point outside \(Y\), the preceding rank difference
is multiplied by a fibre matrix size of \(I_{r,s}\), namely
\(r,s\), or \(rs\). Consequently,
\[
 \rank(b_1\otimes1)(x)-\rank(a_2\otimes1)(x)
 \ge\min\{r,s\}>d+1
 \qquad(x\in X\setminus Y).
\]

Apply \cite[Corollary~3.3]{RT} successively through the
pullback decomposition of \(\widetilde R\otimes I_{r,s}\),
as in the proof of \cite[Corollary~3.4]{RT}, preserving
the values prescribed by \(z_0\).
At each stage, apply the extension result relative to
the union of \(Y\) and the attaching subset. The prescribed
sections agree on their intersection because the preceding
stage agrees with the corresponding image of \(z_0\)
in the quotient. We obtain
\(z\in(\widetilde R\otimes I_{r,s})\otimes\K\) such that
\[
 z^*z=a_2\otimes1,\qquad
 zz^*\in\operatorname{Her}(b_1\otimes1).
\]
\(a_2\otimes1_{I_{r,s}}\precsim b_1\otimes1_{I_{r,s}}\).

Using the unital embedding \(I_{r,s}\hookrightarrow\Zs\)
again, we obtain
\[
 (a-\varepsilon)_+\otimes1_{\Zs}
 \precsim a_2\otimes1_{\Zs}
 \precsim b_1\otimes1_{\Zs}
 \precsim b\otimes1_{\Zs}.
\]
If a unit was adjoined, comparison between these elements
is unchanged upon restricting to the ideal \(A\otimes\K\otimes\Zs\).
Since \(\varepsilon>0\) was arbitrary,
\(a\otimes1_{\Zs}\precsim b\otimes1_{\Zs}\).
\end{proof}
\begin{proof}[Proof of Theorem~\ref{thm:D}]
Let \(x=[a]\) and \(y=[b]\) in \(\Cu(A)\) satisfy
\[
 2x\le 2y,\qquad 3x\le3y
\]
for some \(k\ge1\). Proposition~\ref{prop:Ztest} gives
\(a\otimes1_{\Zs}\precsim b\otimes1_{\Zs}\).
Since \(A\) is separable and \(\Zs\)-stable, the first-factor
inclusion \(A\to A\otimes\Zs\) is approximately unitarily equivalent
to an isomorphism \cite[Theorem~2.2]{TW}. Hence \(a\precsim b\),
so \(x\le y\). The characterization of near unperforation recalled
in Section~2 proves the theorem.
\end{proof}
\section*{ AI statement}
During the preparation of this manuscript, the authors used ChatGPT 6 Astra to assist with checking proofs, locating references, and revising the text. All mathematical arguments, references, and edits were independently reviewed and verified by the authors, who take full responsibility for the content of the paper. 


\begingroup\small
\begin{thebibliography}{99}

\bibitem{AEGL}
Q.~An, S.~Eilers, G.~Gong and Z.~Liu, \doi{10.1007/s00208-026-03603-5}{Subhomogeneity in the classification of real rank zero \(C^*\)-algebras}, \emph{Math. Ann.} \textbf{396} (2026), no.~3, article~71.


\bibitem{APT}
R. Antoine, F. Perera and H. Thiel,
\doi{10.1090/memo/1199}{Tensor products and regularity properties of Cuntz semigroups},
\emph{Mem. Amer. Math. Soc.} \textbf{251} (2018), no.~1199.

\bibitem{Bla}
B. Blackadar,
\href{https://bruceblackadar.com/Mathematics/book6.pdf}{\emph{\(K\)-Theory for Operator Algebras}},
2nd ed., Math. Sci. Res. Inst. Publ., vol.~5, Cambridge University Press, 1998.

\bibitem{BR92}
B. Blackadar and M. R{\o}rdam,
\doi{10.1016/0021-8693(92)90098-7}{Extending states on preordered semigroups and the existence of quasitraces on \(C^*\)-algebras},
\emph{J. Algebra} \textbf{152} (1992), 240--247.

\bibitem{Br}
L.~G. Brown, \doi{10.2140/pjm.1977.71.335}{Stable isomorphism of hereditary subalgebras of \(C^*\)-algebras}, \emph{Pacific J. Math.} \textbf{71} (1977), 335--348.

\bibitem{BrownAF}
L.~G. Brown, Extensions of AF algebras: the projection lifting problem,
in \emph{Operator Algebras and Applications, Part~1}
(R.~V. Kadison, ed.), Proc. Sympos. Pure Math., vol.~38,
Amer. Math. Soc., Providence, RI, 1982, pp.~175--176.

\bibitem{BDF}
L.~G. Brown, R.~G. Douglas and P.~A. Fillmore,
\doi{10.2307/1970999}{Extensions of \(C^*\)-algebras and \(K\)-homology},
\emph{Ann. of Math.} (2) \textbf{105} (1977), 265--324.

\bibitem{LBD}
L.~G. Brown and M. Dadarlat,
\doi{10.1112/jlms/53.3.582}{Extensions of \(C^*\)-algebras and quasidiagonality},
\emph{J. London Math. Soc.} (2) \textbf{53} (1996), 582--600.

\bibitem{BP}
L.~G. Brown and G.~K. Pedersen, \doi{10.1016/0022-1236(91)90056-B}{\(C^*\)-algebras of real rank zero}, \emph{J. Funct. Anal.} \textbf{99} (1991), 131--149.

\bibitem{BrownQD}
N.~P. Brown, \doi{10.2969/aspm/03810019}{On quasidiagonal \(C^*\)-algebras},
in \emph{Operator Algebras and Applications}, Adv. Stud. Pure Math., vol.~38,
Math. Soc. Japan, Tokyo, 2004, pp.~19--64.

\bibitem{NBD}
N.~P. Brown and M. Dadarlat,
\doi{10.2969/aspm/03810065}{Extensions of quasidiagonal \(C^*\)-algebras and \(K\)-theory},
in \emph{Operator Algebras and Applications}, Adv. Stud. Pure Math., vol.~38,
Math. Soc. Japan, 2004, pp.~65--84.

\bibitem{Busby}
R.~C. Busby, \doi{10.2307/1994883}{Double centralizers and extensions of \(C^*\)-algebras},
\emph{Trans. Amer. Math. Soc.} \textbf{132} (1968), 79--99.

\bibitem{CEI}
K.~T. Coward, G.~A. Elliott and C.~Ivanescu,
\doi{10.1515/CRELLE.2008.075}{The Cuntz semigroup as an invariant for \(C^*\)-algebras},
\emph{J. Reine Angew. Math.} \textbf{623} (2008), 161--193.

\bibitem{Cun78}
J. Cuntz, \doi{10.1007/BF01421922}{Dimension functions on simple \(C^*\)-algebras},
\emph{Math. Ann.} \textbf{233} (1978), 145--153.

\bibitem{DadUCT}
 M. Dadarlat,
 Some remarks on the universal coefficient theorem in \(KK\)-theory,
 in \emph{Operator Algebras and Mathematical Physics}
 (Constan\c{t}a, 2001), Theta, Bucharest, 2003, pp.~65--74.

\bibitem{Eck}
C. Eckhardt, \href{https://arxiv.org/abs/2611.14888}{Strongly quasidiagonal \(C^*\)-algebras with nuclear dimension one and infinite decomposition rank}, arXiv:2611.14888 (2026).

\bibitem{EHS}
E.~G. Effros, D.~E. Handelman and C.-L. Shen,
\doi{10.2307/2374244}{Dimension groups and their affine representations},
\emph{Amer. J. Math.} \textbf{102} (1980), 385--407.

\bibitem{ELP}
S.~Eilers, T.~A. Loring and G.~K. Pedersen, \doi{10.1515/crll.1998.055}{Stability of anticommutation relations. An application of noncommutative CW complexes}, \emph{J. Reine Angew. Math.} \textbf{499} (1998), 101--143.

\bibitem{ERR}
S. Eilers, G. Restorff and E. Ruiz,
\doi{10.4153/CJM-2013-015-7}{The ordered \(K\)-theory of a full extension},
\emph{Canad. J. Math.} \textbf{66} (2014), 596--624.

\bibitem{Ell76}
G.~A. Elliott,
\doi{10.1016/0022-1236(76)90054-9}{Automorphisms determined by multipliers on ideals of a \(C^*\)-algebra},
\emph{J. Funct. Anal.} \textbf{23} (1976), 1--10.

\bibitem{Ell93}
G.~A. Elliott, \doi{10.1515/crll.1993.443.179}{On the classification of \(C^*\)-algebras of real rank zero},
\emph{J. Reine Angew. Math.} \textbf{443} (1993), 179--219.

\bibitem{EG}
G.~A. Elliott and G. Gong, \doi{10.2307/2118565}{On the classification of \(C^*\)-algebras of real rank zero, II}, \emph{Ann. of Math.} (2) \textbf{144} (1996), 497--610.

\bibitem{GR}
J. Gabe and E. Ruiz, \doi{10.1017/S0017089519000053}{The unital Ext-groups and classification of \(C^*\)-algebras}, \emph{Glasgow Math. J.} \textbf{62} (2020), 201--231.

\bibitem{Goo}
K.~R. Goodearl, Notes on a class of simple \(C^*\)-algebras with real rank zero, \emph{Publ. Mat.} \textbf{36} (1992), no.~2A, 637--654.

\bibitem{Had}
D.~Hadwin, \href{https://jot.theta.ro/jot/archive/1987-018-001/index_1987-018-001.html}{Strongly quasidiagonal \(C^*\)-algebras}, \emph{J. Operator Theory} \textbf{18} (1987), 3--18.

\bibitem{JS99}
X. Jiang and H. Su,
\doi{10.1353/ajm.1999.0012}{On a simple unital projectionless \(C^*\)-algebra},
\emph{Amer. J. Math.} \textbf{121} (1999), 359--413.

\bibitem{LR}
H. Lin and M. R{\o}rdam,
\doi{10.1112/jlms/51.3.603}{Extensions of inductive limits of circle algebras},
\emph{J. London Math. Soc.} (2) \textbf{51} (1995), 603--613.

\bibitem{NW}
P.~W. Ng and W. Winter, \href{https://arxiv.org/abs/math/0601069}{A note on subhomogeneous \(C^*\)-algebras}, \emph{C.~R. Math. Rep. Acad. Sci. Canada} \textbf{28} (2006), 91--96.

\bibitem{ORT}
E. Ortega, M. R{\o}rdam and H. Thiel,
\doi{10.1016/j.jfa.2011.02.017}{The Cuntz semigroup and comparison of open projections},
\emph{J. Funct. Anal.} \textbf{260} (2011), 3474--3493.

\bibitem{PRRS}
D. Pask, I. Raeburn, M. R{\o}rdam and A. Sims, \href{https://arxiv.org/abs/math/0512254}{Rank-two graphs whose \(C^*\)-algebras are direct limits of circle algebras}, \emph{J. Funct. Anal.} \textbf{239} (2006), 137--178.

\bibitem{PhillipsLim}
N.~C. Phillips, \doi{10.1090/S0002-9947-07-03849-4}{Cancellation and stable rank for direct limits of recursive subhomogeneous algebras}, \emph{Trans. Amer. Math. Soc.} \textbf{359} (2007), 4625--4652.

\bibitem{Phillips}
N.~C. Phillips, \href{https://arxiv.org/abs/math/0101156}{Recursive subhomogeneous algebras}, \emph{Trans. Amer. Math. Soc.} \textbf{359} (2007), 4595--4623.

\bibitem{Rie}
M.~A. Rieffel,
\doi{10.1112/plms/s3-46.2.301}{Dimension and stable rank in the \(K\)-theory of \(C^*\)-algebras},
\emph{Proc. London Math. Soc.} (3) \textbf{46} (1983), 301--333.

\bibitem{RT}
L. Robert and A. Tikuisis,
\doi{10.1112/plms/pdq017}{Hilbert \(C^*\)-modules over a commutative \(C^*\)-algebra},
\emph{Proc. Lond. Math. Soc.} (3) \textbf{102} (2011), 229--256.

\bibitem{Ror04}
M. R{\o}rdam,
\doi{10.1142/S0129167X04002661}{The stable and the real rank of \(\Zs\)-absorbing \(C^*\)-algebras},
\emph{Internat. J. Math.} \textbf{15} (2004), 1065--1084.

\bibitem{RW}
M. R{\o}rdam and W. Winter,
\doi{10.1515/crelle.2010.039}{The Jiang--Su algebra revisited},
\emph{J. Reine Angew. Math.} \textbf{642} (2010), 129--155.

\bibitem{Sp}
J.~S. Spielberg,
\doi{10.1016/0022-1236(88)90104-8}{Embedding \(C^*\)-algebra extensions into AF algebras},
\emph{J. Funct. Anal.} \textbf{81} (1988), 325--344.

\bibitem{TW}
A.~S. Toms and W. Winter,
\doi{10.1090/S0002-9947-07-04173-6}{Strongly self-absorbing \(C^*\)-algebras},
\emph{Trans. Amer. Math. Soc.} \textbf{359} (2007), 3999--4029.

\bibitem{TWASH}
A.~S. Toms and W. Winter, \doi{10.4153/CJM-2008-031-6}{\(\Zs\)-stable ASH algebras},
\emph{Canad. J. Math.} \textbf{60} (2008), 703--720.

\end{thebibliography}
\endgroup
\end{document}